\documentclass[11pt, doublespaced, leqno,noinfoline]{article}

\RequirePackage[OT1]{fontenc}

\RequirePackage[aop]{imsart}

\usepackage{amsmath,amsthm,amssymb,graphicx,natbib,color}

\startlocaldefs
\def\@journal{\ }

\newtheorem{thm}{Theorem}[section]
\newtheorem{remark}[thm]{Remark}
\newtheorem{propn}[thm]{Proposition}

\newtheorem{cor}[thm]{Corollary}

\numberwithin{equation}{section}

\endlocaldefs

\begin{document}

\begin{frontmatter}


\title{MAHARAM EXTENSION FOR NONSINGULAR GROUP ACTIONS}

\runtitle{Maharam Extension for Actions}


\author{\fnms{Parthanil} \snm{Roy}\ead[label=e1]{pr72@cornell.edu}}
\address{School of Operations Research\\ and Industrial Engineering\\Cornell University\\Ithaca, NY 14853\\\printead{e1}}\and
\author{\fnms{Gennady} \snm{Samorodnitsky}\ead[label=e2]{gennady@orie.cornell.edu}}
\address{School of Operations Research\\ and Industrial Engineering\\Cornell University\\Ithaca, NY 14853\\\printead{e2}}\thanksref{t1}
\thankstext{t1}{Supported in part by NSF
grant DMS-0303493, NSA grant MSPF-05G-049
 and NSF training grant ``Graduate and Postdoctoral
Training in Probability and Its Applications'' at Cornell
University.}

\affiliation{Cornell University}

\runauthor{P.~Roy and G.~Samorodnitsky}

\begin{abstract}

We establish a generalization of the Maharam Extension Theorem to
nonsingular $\mathbb{Z}^d$-actions. We also present an extension of
Krengel's representation of dissipative transformations to nonsingular
actions.

\end{abstract}



\end{frontmatter}

\section{Introduction}

Maharam extension theorem extends in a natural way an invertible
nonsingular conservative transformation of a $\sigma$-finite
standard measure space to an invertible conservative {\it measure
preserving} transformation on an extended space, the so-called {\it
Maharam skew product}.  The result was established in
\cite{maharam:1964}) and has been used in a number of ways,
allowing, in particular, extensions of certain notions from the
measure preserving case to the nonsingular case; see e.g.
\cite{silva:thieullen:1995}.

In this paper we generalize Maharam's theorem to nonsingular
$\mathbb{Z}^d$-actions. Our approach  is  different than the one
often used in the case $d=1$, based on the fact that conservativity
is equivalent to incompressibility. We use, instead, a result on the
maximal value assigned by a group action over an increasing sequence
of cubes to a nonnegatve function (Proposition \ref{pr:max} below),
which may be of an independent interest.  In the proof of one of
the statements in that proposition we use a recently established
extension of Krengel's theorem (see \cite{krengel:1969}) on the
structure of dissipative nonsingular transformations to nonsingular
$\mathbb{Z}^d$-actions. This result has not, apparently, been stated
before. Apart from that, the proof of the main result of this paper
is entirely from the first principles.

We state the extensions of both Maharam's theorem and Krengel's
theorem in Section \ref{sec_prelim}. The proof of Maharam Extension
Theorem is given in Section $\ref{sec_proof}$.

\section{Maharam's theorem and Krengel's theorem for nonsingular
$\mathbb{Z}^d$-actions } \label{sec_prelim}

Let $\{\phi_t\}_{t \in \mathbb{Z}^d}$ be a nonsingular action on a
standard Borel space $(S, \mathcal{S})$ with a $\sigma$-finite
measure $\mu$. Then, by Theorem $1$ in \cite{maharam:1964},
\[
\phi^{\ast}_t(s,y):=(\phi_t(s), y\frac{d\mu}{d\mu\circ
\phi_t}(s)),\;\;\;\; t \in \mathbb{Z}^d
\]
is a measure preserving group action on the product space $\bigl( S
\times (0,\infty), \mathcal{S}\times \mathcal{B},
\mu\times\text{Leb}\bigr)$. Here Leb is the Lebesgue measure on
$(0,\infty)$.

The following is our main result.
\begin{thm}\label{thm_maharam_extn}
The group action   $\{\phi^{\ast}_t\}_{t \in
  \mathbb{Z}^d}$ is conservative on $\bigl( S \times (0,\infty),
  \mathcal{S}\times \mathcal{B}, \mu\times\text{Leb}\bigr)$ if and
  only if the group action  $\{\phi_t\}_{t \in
  \mathbb{Z}^d}$ is conservative on $(S, \mathcal{S}, \mu)$.
\end{thm}
In the case $d=1$ this is the content of Maharam Extension Theorem
(\cite{maharam:1964}).

The proof of Theorem \ref{thm_maharam_extn} presented in the next
section relies on a result on the maximal value of a function
transformed by the group of dual operators, given in Proposition
\ref{pr:max}. The argument for one part of the proposition uses the
following extension of Krengel's Theorem (see \cite{krengel:1969})
on the structure of dissipative nonsingular maps to
$\mathbb{Z}^d$-actions. It follows immediately from Theorem $2.2$ in
\cite{rosinski:2000} and Corollary $2.4$ in
\cite{roy:samorodnitsky:2006}. It appears that the result has not
been stated previously. Recall that nonsingular group actions
$\{\phi_t\}_{t\in G}$ and $\{\psi_t\}_{t\in G}$, defined on standard
measure spaces $(S,\mathcal{S},\mu)$ and $(T,\mathcal{T},\nu)$
resp., are equivalent if there is a Borel isomorphism $\Phi$ between
the measure spaces such that $\nu\sim \mu \circ \Phi^{-1}$ and for
each $t \in T$, $\psi_t \circ \Phi=\Phi \circ \phi_t \;$  $\mu$-a.e.
\begin{thm}[Krengel's Theorem for
    $\mathbb{Z}^d$-actions]\label{thm_krengel_Z^d}
Let $\{\phi_t\}$ be a nonsingular $\mathbb{Z}^d$-action on a
$\sigma$-finite standard measure space $(S,\mathcal{S},\mu)$. Then
$\{\phi_t\}$ is dissipative if and only if it is equivalent to the
$\mathbb{Z}^d$-action
\begin{equation}
\psi_t(w,s):=(w,t+s),\;\;\;t\in \mathbb{Z}^d
\label{form_of_diss_action}
\end{equation}
defined on $(W \times {\mathbb{Z}}^d, \tau \otimes l)$, where $(W,
\mathcal{W},\tau)$ is some $\sigma$-finite standard measure space
and $l$ is the counting measure on $\mathbb{Z}^d$.
\end{thm}

\section{Proof of Theorem \ref{thm_maharam_extn}} \label{sec_proof}

Let $\{\phi_t\}_{t \in \mathbb{Z}^d}$ be as above and
$\hat{\phi}_t:\,L^1(\mu)\rightarrow L^1(\mu)$ be the dual to
$\phi_{-t}$ operator (see Section $1.3$ in \cite{aaronson:1997})
\[
\hat{\phi}_t g(s)=g\circ\phi_t(s)\,\frac{d\mu \circ \phi_t}{d
\mu}(s)\,,\;\;\;\;s \in S\,.
\]
The following result, which may be of independent interest, is the
key step in the proof of Theorem $\ref{thm_maharam_extn}$. The
inequalities in the statement of this proposition and elsewhere are
understood in the sense of the natural partial order on
$\mathbb{Z}^d$.
\begin{propn}\label{pr:max}
\noindent (a) If $\{\phi_t\}_{t \in \mathbb{Z}^d}$ is conservative
then for all $g \in L^1(\mu),\,g \geq 0$, we have
\begin{equation}
\frac{1}{n^d} \int_S \max_{0 \leq t \leq (n-1)\mathbf{1}}
\hat{\phi}_t g(s) \mu(ds) \rightarrow 0\,.\label{max_cons}
\end{equation}

\noindent (b) If $\{\phi_t\}_{t \in \mathbb{Z}^d}$ is dissipative
then for all $g \in L^1(\mu),\,g \geq 0,\,\mu(g>0)>0$, we have
\begin{equation}
\frac{1}{n^d} \int_S \max_{0 \leq t \leq (n-1)\mathbf{1}}
\hat{\phi}_t g(s) \mu(ds) \rightarrow a \label{max_diss}
\end{equation}
for some $0<a<\infty$.
\end{propn}

\begin{proof}
(a) There is no loss of generality in assuming that $\mu$ is a
probability measure. We can also assume that the support of the
family $\bigl\{\hat{\phi}_t g\bigr\}_{t \in \mathbb{Z}^d}$ is the
entire set $S$. Let $\{\alpha_u:\,u \in \mathbb{Z}^d\}$ be a
collection of positive numbers summing up to $1$. Then applying the
group action version of Theorem 1.6.3 in \cite{aaronson:1997} to
$f=\sum_{u \in \mathbb{Z}^d}\alpha_u\,\hat{\phi}_u g$ we have,
\begin{equation}
\sum_{t \in \mathbb{Z}^d}\hat{\phi}_t g(s)=\infty \;\;\;\;\mbox{for
}\mu\mbox{-a.a.}\;s\,.\label{eqn_sum_inf}
\end{equation}

To prove $(\ref{max_cons})$ we will show that
\[
a_n:=\frac{1}{(2n+1)^d} \int_S \max_{t \in J_n} \hat{\phi}_t g(s)
\mu(ds) \rightarrow 0\,,
\]
where, $J_n:=\{(i_1,i_2,\ldots,i_d):-n \leq i_1,i_2,\ldots,i_d \leq
n \}$. Note that
\begin{eqnarray*}
a_n&\leq&\frac{1}{(2n+1)^d}\bigg(\int_S \max_{t \in J_n}
\Big[\hat{\phi}_t g(s)I\Big(\hat{\phi}_t g(s)\leq \epsilon \sum_{u
\in
J_n} \hat{\phi}_u g(s)\Big)\Big] \mu(ds)\\
&&\;\;\;\;\;\;\;\;\;\;\;\;\;\;\;\;\;+\int_S \max_{t \in J_n}
\Big[\hat{\phi}_t g(s)I\Big(\hat{\phi}_t g(s)> \epsilon \sum_{u \in
J_n} \hat{\phi}_u g(s)\Big)\Big] \mu(ds)\bigg)\\
&&=a^{(1)}_n+a^{(2)}_n\,,
\end{eqnarray*}
where $\epsilon > 0$ is arbitrary. Clearly,
\begin{equation}
a^{(1)}_n \leq \frac{\epsilon}{(2n+1)^d} \sum_{u \in J_n} \int_S
\hat{\phi}_u g(s) \mu(ds) =\epsilon \|g\| \;,
\label{bound_on_a^{(1)}_n}
\end{equation}
where, $\|g\|:=\int_S g(s)\mu(ds)<\infty$\,. Also, by duality,
\begin{eqnarray}
a^{(2)}_n &\leq&  \frac{1}{(2n+1)^d} \sum_{t \in J_n} \int_S
\hat{\phi}_t g(s)I_{A_{t,n}}(s) \mu(ds)
\label{bound_on_a^{(2)}_n} \\
&=& \frac{1}{(2n+1)^d} \sum_{t \in J_n}\int_S g(s)
I_{\phi_t^{-1}(A_{t,n})}(s) \mu(ds)\,,\nonumber
\end{eqnarray}
where, $A_{t,n}=\{s:\, \hat{\phi}_t g(s) > \epsilon \sum_{u \in J_n}
\hat{\phi}_u g(s)\}\,,\,n \geq 1,\, t \in J_n\,$. Define
\[
U_n:=\{(i_1,i_2,\ldots,i_d): -n+[\sqrt{n}] \leq i_1,i_2,\ldots,i_d
\leq n-[\sqrt{n}]\} \,.
\]
Observe that by the nonnegativity, for every $t \in U_n$,
\begin{eqnarray*}
\phi_t^{-1}(A_{t,n}) &=&\Big\{s:\,g(s)>\epsilon \sum_{u \in J_n}
\hat{\phi}_{u+t} g(s) \Big\}\\
\;\;\;\;&\subseteq& \Big\{s:\,g(s)>\epsilon \sum_{u \in
J_{[\sqrt{n}]}} \hat{\phi}_{u} g(s)\Big\} \,.
\end{eqnarray*}
Therefore, for any $M>0$
$$
\max_{t \in U_n}\;\mu(\phi_t^{-1}(A_{t,n}))\leq
\mu\{s:\,g(s)>\epsilon M\}+\mu\bigg(\sum_{t \in J_{[\sqrt{n}]}}
\hat{\phi}_t g(s)\leq M\bigg) \,.
$$
Letting first $n\to\infty$, using $(\ref{eqn_sum_inf})$, and then letting
$M\to\infty$ we see that
$$
\lim_{n \rightarrow \infty} \max_{t \in
U_n}\;\mu(\phi_t^{-1}(A_{t,n}))=0\;. \label{limit_of_max_mu}
$$
From here we immediately see that
\begin{eqnarray}
&&\frac{1}{(2n+1)^d} \sum_{t \in U_n} \int_S \hat{\phi}_t g(s)I_{A_{t,n}}(s) \mu(ds) \nonumber \\
&&\;\;\;\;\;\;\;\;\;\;\;\;\;\;\;\;\;\;\;\;\;\;\;\;\;\;\;\;=\frac{1}{(2n+1)^d}
\sum_{t \in U_n}\int_{\phi_t^{-1}(A_{t,n})} g(s) \mu(ds) \rightarrow
0\,. \label{bound_on_part_of_a^{(2)}_n}
\end{eqnarray}
Define $V_n=J_n\setminus U_n$, and note that Card$(V_n)=o(n^d)$ as
$n\to\infty$.  Therefore, using $(\ref{bound_on_a^{(2)}_n})$ and
$(\ref{bound_on_part_of_a^{(2)}_n})$ we have,
$$
a^{(2)}_n \leq \frac{1}{(2n+1)^d} \sum_{t \in U_n} \int_S \hat{\phi}_t
g(s)I_{A_{t,n}}(s) \mu(ds)+\,\frac{\text{Card$(V_n)$}}{(2n+1)^d}\,
\|g\| \rightarrow 0\,,
$$
implying that
\[
\limsup a_n \leq \limsup a^{(1)}_n + \limsup a^{(2)}_n \leq \epsilon
\|g\|\,.
\]
Since $\epsilon >0$ is arbitrary, the claim follows.\\

\noindent (b) Since the statement is invariant under a passage from
one group action to an equivalent one, we will use Theorem
\ref{thm_krengel_Z^d} and check that for any $\sigma$-finite standard
measure space $(W, \mathcal{W},\tau)$ we have for all $f \in L^1(W
\times \mathbb{Z}^d, \tau \otimes l)$ and $f\geq 0$ with $\tau \otimes
l(f>0)>0$,
\begin{equation}
\frac{1}{n^d}\sum_{s \in \mathbb{Z}^d} \int_W \max_{0 \leq t
\leq(n-1)\mathbf{1}} f(w,s+t)\, \tau(dw) \rightarrow a
\label{max_diss_krengel}
\end{equation}
for some $0<a<\infty$. In fact, we will show that
$\eqref{max_diss_krengel}$ holds with $a=\int_W h(w)\tau(dw) \in
(0,\infty)$ where $h(w):=\sup_{s\in \mathbb{Z}^d}f(w,s)$ for all $w
\in W$.

We start with the case where
$f$ has compact support, that is
\[
f(w,s)I_{W \times [-m\mathbf{1},m\mathbf{1}]^{c}}(w,s) \equiv 0
\mbox{ for some }m=1,2,\ldots,
\]
where $[u,v]:=\{t \in \mathbb{Z}^d: u \leq t \leq v\}$.
In that case, we have, for all $n \geq 2m-1$,
\begin{eqnarray*}
&&\sum_{s \in \mathbb{Z}^d} \int_W \max_{0 \leq t
\leq(n-1)\mathbf{1}} f(w,s+t)\,\tau(dw)\\
            &=& \sum_{(-m-n+1)\mathbf{1} \leq s \leq m\mathbf{1}} \int_W
\max_{0 \leq t \leq(n-1)\mathbf{1}}f(w,s+t)\,\tau(dw)\\
            &=& \sum_{s \in A_n} \int_W \max_{0 \leq t
\leq(n-1)\mathbf{1}}f(w,s+t)\,\tau(dw)\\
            &&+ \sum_{s \in B_n} \int_W \max_{0 \leq t
\leq(n-1)\mathbf{1}}f(w,s+t)\,\tau(dw)=:T_n + R_n\,,
\end{eqnarray*}
where $A_n=[(m-n-1)\mathbf{1},-m\mathbf{1}]$  and
$B_n=[(-m-n+1)\mathbf{1},m\mathbf{1}]-[(m-n-1)\mathbf{1},-m\mathbf{1}]$.
Observe that, for $n \geq 2m+1$ we have for each $s \in
A_n$,
\[
\max_{0 \leq t \leq(n-1)\mathbf{1}}f(w,s+t)|=h(w)
\]
while for each $s \in B_n$
\[
\max_{0 \leq t \leq(n-1)\mathbf{1}}f(w,s+t)| \leq h(w)\,,
\]
and so
\begin{eqnarray*}
T_n&=&(n-2m)^d \int_W h(w) \tau(dw)\,,\\
R_n&\leq& [(2m+n)^d-(n-2m)^d] \int_W h(w) \tau(dw).
\end{eqnarray*}
Therefore $(\ref{max_diss_krengel})$ follows when $f$ has compact
support. In the general case, given $\epsilon > 0$, choose a compactly
supported $f_{\epsilon}$ such that $f_{\epsilon}(w,s) \leq f(w,s)$
for all $w,s$ and
\[
\sum_{s \in \mathbb{Z}^d} \int_W f(w,s)\,\tau(dw)-\sum_{s \in
\mathbb{Z}^d} \int_W f_{\epsilon}(w,s)\,\tau(dw) \leq \epsilon .
\]
Let
\[
h_\epsilon (w)= \sup_{s \in \mathbb{Z}^d} f_{\epsilon}(w,s),\;\;\;w
\in W.
\]
Then
\begin{eqnarray*}
0 &\leq& \int_W h(w) \tau(dw)-\int_W h_\epsilon(w) \tau(dw)
\\
 &\leq& \int_W \sup_{s \in \mathbb{Z}^d}
\bigl(f(w,s)-f_\epsilon(w,s) \bigr)\,\tau(dw)\\
 &\leq& \int_W \sum_{s \in \mathbb{Z}^d}
\bigl(f(w,s)-f_\epsilon(w,s) \bigr)\,\tau(dw)\\
 &=& \sum_{s \in \mathbb{Z}^d} \int_W f(w,s)\,\tau(dw)-\sum_{s \in
\mathbb{Z}^d} \int_W f_{\epsilon}(w,s)\,\tau(dw) \leq \epsilon.
\end{eqnarray*}
Therefore,
\begin{eqnarray*}
&&\Bigl|\frac{1}{n^d}\sum_{s \in \mathbb{Z}^d} \int_W
\max_{0 \leq t \leq(n-1)\mathbf{1}}f(w,s+t)\,\tau(dw)-\int_W h(w)\,\tau(dw)\Bigr|\\
                         &\leq& \frac{1}{n^d}\Bigl|\sum_{s \in \mathbb{Z}^d} \int_W
\max_{0 \leq t \leq(n-1)\mathbf{1}}f(w,s+t)\,\tau(dw)\\
                         &&\;\;\;\;\;\;\;\;\;\;-\sum_{s \in \mathbb{Z}^d}
\int_W \max_{0 \leq t \leq(n-1)\mathbf{1}}f_\epsilon (w,s+t)\,
\tau(dw)\Bigr|\\
                         &+&\Bigl|\frac{1}{n^d}\sum_{s \in \mathbb{Z}^d} \int_W
\max_{0 \leq t \leq(n-1)\mathbf{1}}f_\epsilon (w,s+t)\,
\tau(dw)-\int_W
h_\epsilon(w) \tau(dw)\Bigr|\\
                         &+&\Bigl|\int_W h_\epsilon(w) \tau(dw)-\int_W
h(w) \tau(dw)\Bigr|=:T_n^{(1)}+T_n^{(2)}+T_n^{(3)}.
\end{eqnarray*}
By the above, $T_n^{(3)}\leq \epsilon$, and the same argument shows
that $T_n^{(1)}\leq \epsilon$ as well. Furthermore, by the already
considered compact support case, $T_n^{(2)} \rightarrow 0$ as $n
\rightarrow \infty$. Hence
\[
\limsup_{n\rightarrow \infty} |\frac{1}{n^d}\sum_{s \in
\mathbb{Z}^d} \int_W \max_{0 \leq t
\leq(n-1)\mathbf{1}}f(w,s+t)\,\tau(dw)-\int_W h(w)\,\tau(dw)| \leq
2\epsilon,
\]
and, since $\epsilon > 0$ is arbitrary, the proof is complete.

\end{proof}

The following corollary is immediate.
\begin{cor} \label{cor_noncons}
If $g \in L^1(\mu),\,g \geq
0$, and $\mu(\mbox{Support}(g) \cap \mathcal{D})
>0$ where $\mathcal{D}$ is the dissipative part of $\{\phi_t\}$, then
\[
\frac{1}{n^d} \int_S \max_{0 \leq t \leq (n-1)\mathbf{1}}
\hat{\phi}_t g(s) \mu(ds) \rightarrow a
\]
for some $0<a<\infty$.
\end{cor}

\begin{remark}\label{remark_cons}
\textnormal{From Corollary $\ref{cor_noncons}$ it follows that, if
$(\ref{max_cons})$ holds for some $g \in L^1(\mu),\,g \geq 0$, then
\[
\mbox{Support}(g) \subseteq \mathcal{C}\;\;\;\mbox{mod}\;\mu\,,
\]
where $\mathcal{C}$ is the conservative part of $\{\phi_t\}$. In
other words, if there exists a sequence of functions $g_m \in
L^1(\mu),\,g_m \geq 0$, whose support increases to $S$,
such that $(\ref{max_cons})$ holds for $g_m$ for all $m \geq 1$,
then $\{\phi_t\}$ is conservative.}
\end{remark}

\begin{proof}[Proof of Theorem $\ref{thm_maharam_extn}$] If
  $\{\phi^{\ast}_t\}$ is conservative, so is clearly
  $\{\phi_t\}$. Suppose now that $\{\phi_t\}$ is conservative. To show
  conservativity of $\{\phi^{\ast}_t\}$
we will use Remark $\ref{remark_cons}$. Since
$\mu$ is $\sigma$-finite, there is a sequence of measurable sets
$S_m \uparrow S$, such that, $\mu(S_m) <\infty$ for all $m \geq 1$.
Consider a sequence of nonnegative functions $g^{\ast}_m:=I_{S_m
  \times (0,m)} \in L^1(\mu \otimes Leb)$, $m\geq 1$. Note that the
  support of $g_m^\ast$ is $S_m \times (0, m) \uparrow S
\times (0,\infty)\,.$

Observe that $g^{\ast}_m(s,y)=I\{(s,y):\,0
<y<mI_{S_m}(s)\}$. If $w_t:=\frac{d\mu \circ \phi_t}{d\mu},\;t \in
\mathbb{Z}^d$, then for all $m \geq 1$ we have,
\begin{align*}
&\frac{1}{n^d}\int_0^\infty\int_S \max_{0 \leq t \leq
(n-1)\mathbf{1}}
\hat{\phi}^{\ast}_t\,g^{\ast}_m(s,y) \mu(ds)\text{Leb}(dy)\\
&=\frac{1}{n^d}\int_0^\infty\int_S \max_{0 \leq t \leq
(n-1)\mathbf{1}}g^{\ast}_m \circ \phi^{\ast}_t(s,y) \mu(ds)
\text{Leb}(dy)\\
&=\frac{1}{n^d}\int_S \int_0^\infty I\Bigl\{(s,y):\,0<y<\max_{0
    \leq t \leq
    (n-1)\mathbf{1}}mw_t(s)I_{S_m}(\phi_t(s))\Bigr\} \\
&\hskip 3.6in\text{Leb}(dy)\mu(ds)\\
&=\frac{m}{n^d}\int_S \max_{0
\leq t \leq (n-1)\mathbf{1}} \bigg[I_{S_m}(\phi_t(s)) \frac{d\mu
\circ \phi_t}{d\mu}(s)\bigg]
\mu(ds)\\
&=\frac{m}{n^d}\int_S \max_{0 \leq t \leq (n-1)\mathbf{1}}
\hat{\phi}_t I_{S_m}(s) \rightarrow 0
\end{align*}
by part (a) of Proposition \ref{pr:max}. By Remark
$\ref{remark_cons}$ this is enough to prove the theorem.
\end{proof}

\end{document}